\documentclass[conference]{IEEEtran}

\usepackage[noadjust]{cite}
\usepackage{graphicx,color,psfrag,subeqnarray}
\usepackage{amsmath}
\usepackage{amssymb}
\usepackage{amsthm,amssymb}
\usepackage{color}
\usepackage{listings}
\usepackage{multirow}
\usepackage{bm,array}
\usepackage{listings}
\usepackage[font=footnotesize]{caption}
\usepackage{subcaption}
\usepackage{float} 
\usepackage{subfloat}
\usepackage{bigstrut}

\usepackage{algorithm}
\usepackage[noend]{algpseudocode}

\makeatletter
\def\BState{\State\hskip-\ALG@thistlm}
\makeatother

\DeclareCaptionSubType[alph]{figure}
\newcommand\textlcsc[1]{\textsc{\MakeLowercase{#1}}}

\begin{document}

\title{Joint Shaping and Altering the Demand Profile by Residential Plug-in Electric Vehicles for Forward and Spot Markets in Smart Grids}

\author{\IEEEauthorblockN{Farshad Rassaei, Wee-Seng Soh and Kee-Chaing Chua  
		 \\}
\IEEEauthorblockA{Department of Electrical and Computer Engineering\\
National University of Singapore, Singapore\\
Email: farshad@u.nus.edu, \{weeseng,eleckc\}@nus.edu.sg}
}
\maketitle

\begin{abstract}
Plug-in electric vehicles (PEVs) can significantly increase the elasticity of residential electricity demand. This elasticity can be employed to shape the daily aggregated electricity demand profile of a system comprised of a large number of residential PEVs' users sharing one electricity retailer or an aggregator. In this paper, we propose a joint demand shaping and altering algorithm for managing vehicle-to-grid (V2G) enabled PEVs' electricity assignments (charging and discharging) in order to diminish the overall electricity procurement costs for a retailer bidding to two-settlement electricity markets, i.e., a day-ahead (DA) and a real-time (RT) market. This approach is decentralized, scalable, fast converging and does not violate users' privacy. Our simulations' results demonstrate significant overall cost savings (up to 28\%) for a retailer bidding to an operational electricity market by using our proposed algorithm. This becomes even more salient when the power system is integrating a large number of intermittent energy resources wherein RT demand altering is crucial due to more likely contingencies and hence more RT price fluctuations and even more so-called \textit{black swan events}. Lower electricity procurement cost for a retailer finally makes it able to offer better deals to customers and expand its market capacity. This implies that customers can enjoy lower electricity bills as well.   
\end{abstract}


\IEEEpeerreviewmaketitle

\section{Introduction}
\begin{figure}[t]
	\centering
	\includegraphics[width=\columnwidth]{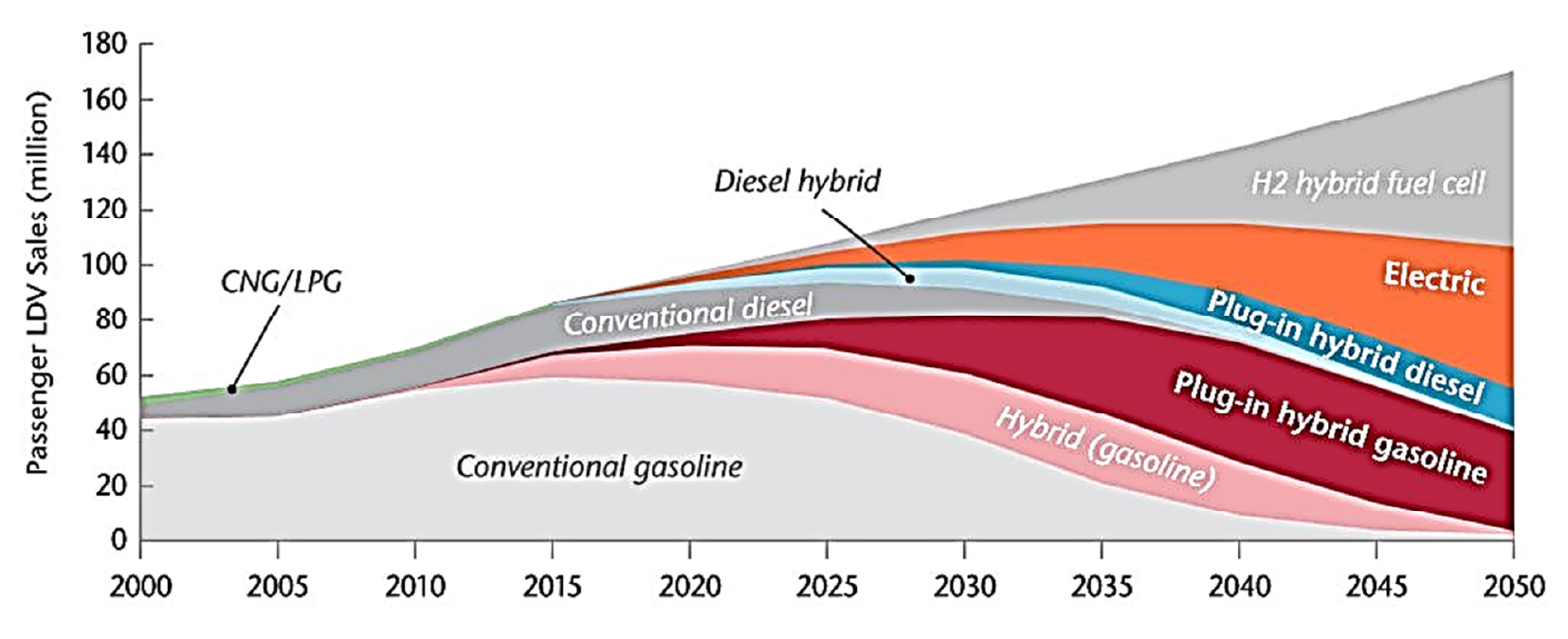}
	\caption{Estimate of Passenger LDV Sales (million) by 2050 \cite{222222}.} \vspace{-1em}
	\label{fevper}
\end{figure}

Ambitious targets and attractive incentives for introducing PEVs into the transport sector have been aimed in many countries \cite{travel2009sustainable}. In fact, the roadmap is that industry and governments attain a combined PEV/PHEV sales share of at least 50\% for light duty vehicle (LDV) sales globally by 2050 \cite{222222}, see Fig. \ref{fevper}.

A generic wholesale electricity market is comprised of electricity markets such as spot and forward markets for trading electricity, and ancillary service markets for guaranteeing security in the provision. In a deregulated electricity system, electricity retailers submit demand bids to wholesale electricity markets. For a day-ahead (DA) electricity market, these demand bids normally have both electricity demand's amount and price components meaning that the retailer buys the specified electricity only if the market clearing price (MCP) is not more than its desired price \cite{mohsenianoptimal}. This bidding can be carried out in a few predefined rounds allowing the retailers to update their bids at each round. This type of bidding is referred to as limit order bidding. Therefore, in this case, the retailer is willing to shape its aggregated demand profile and match it to the electricity profile of a successful bid so that it can minimize its demand from the RT market to balance the load and accordingly reduce the overall electricity procurement cost for each following day.

Based on Pennsylvania–-Jersey-–Maryland (PJM) Interconnection 2014's data \cite{PJM}, the average price of electricity per MWh which has been sold over that year is very close for both DA and RT markets, i.e., \$48.9539 and \$48.2063, respectively. Additionallt, mean reversion theory suggests that prices and returns finally proceed back towards the mean or the average. This mean or average can be the historical average of the price or return or another admissible average \cite{poterba1988mean}. However, hourly pricing data for these DA and RT markets can be significantly distinct and unpredictable at some days and/or hours. Hence, the high uncertainty, particularly in the RT market, can remarkably impact the overall electricity procurement cost for a retailer. This becomes even more salient when the power system is relying on a large number of intermittent energy resources, e.g., wind farms and solar panels and thus more RT price fluctuations and even so-called \textit{black swan} events \cite{taleb2010black} may be occurred.  

In \cite{farisgt15}, we provide a statistical modelling and a closed-form expression for PEVs' uncoordinated charging demand. Furthermore, in \cite{farst1} we propose a decentralized demand shaping algorithm for a priori know desired demand profile for the next day. An overview of DR and their different classifications in a deregulated electricity market is discussed in \cite{albadi2008summary}. The authors in \cite{6910332} compare different bidding rules in wholesale electricity market and in particular when there exist PEVs and renewables' penetration in the power grid. In \cite{mohsenianoptimal}, the author discusses how a time-shiftable load, that may comprise of several smaller time-shiftable subloads, can optimally submit its demand bids to DA and RT markets to minimize its electricity procurement cost. Although this paper provides closed-form solutions for bidding, it does not offer an approach for bids submitted to electricity markets for distributed systems whereby the retailer does not have necessary information about its customers' preferences due to privacy concerns.  

In \cite{kim_bidirectional_2013}, charging and discharging of PEVs are managed in order to maximize the social and individual welfare functions. However, when it comes to residential users, it is not very straightforward to seek pertinent utility and welfare functions for the individual users. Paper \cite{ott2003experience} describes the basic characteristics of the PJM DA and RT electricity markets. The author discusses that economic incentives make the DA and RT market prices converge in the bidding processes. In addition, The locational marginal pricing (LMP) based markets support reliable grid operations through efficient pricing signals to the retailers. Nevertheless, we should note that sending pricing signals and enabling low electricity consumer's bidding in the markets may cause stability problems for the power system, time of use pricing (ToUP) and RT pricing (RTP) methods may lead to high simultaneity in customers' behaviours and instability in the power grid by generating unfavourable aggregated load profiles \cite{roozbehani2012volatility}. In \cite{moghadam2014demand}, a RT pricing mechanism is proposed for the grid operator to manage the DR of retailers or aggregators upon a contingency of supply. 

In \cite{ghamkhari2014optimal}, the authors investigate the bidding process for DA and RT markets for data centres by taking into account both costs and risks. The authors in \cite{6965656} model the bidding interactions between a demand response aggregator (DRA) and generators as a Stackelberg game \cite{von2010market} for a DA market wherein the game is one of complete information. 

In this paper, we propose a model and a fast converging and decentralized joint demand shaping and altering algorithm for managing PEVs' electricity assignments (charging and discharging) in order to minimize the electricity procurement costs of a retailer's bidding to two-settlement electricity markets, i.e., DA and RT markets. We adopt PJM Interconnection \cite{PJM} pricing for 2014 to evaluate our algorithms' results.

\section{System Model} \label{SM}
In this section, we explain the underlying model of the power system in this paper which entails the energy markets, the electricity retailers or the aggregators, and the users. We articulate different parts of this model in the sequel.

\begin{figure}[t]
	\centering
	\includegraphics[width=\linewidth]{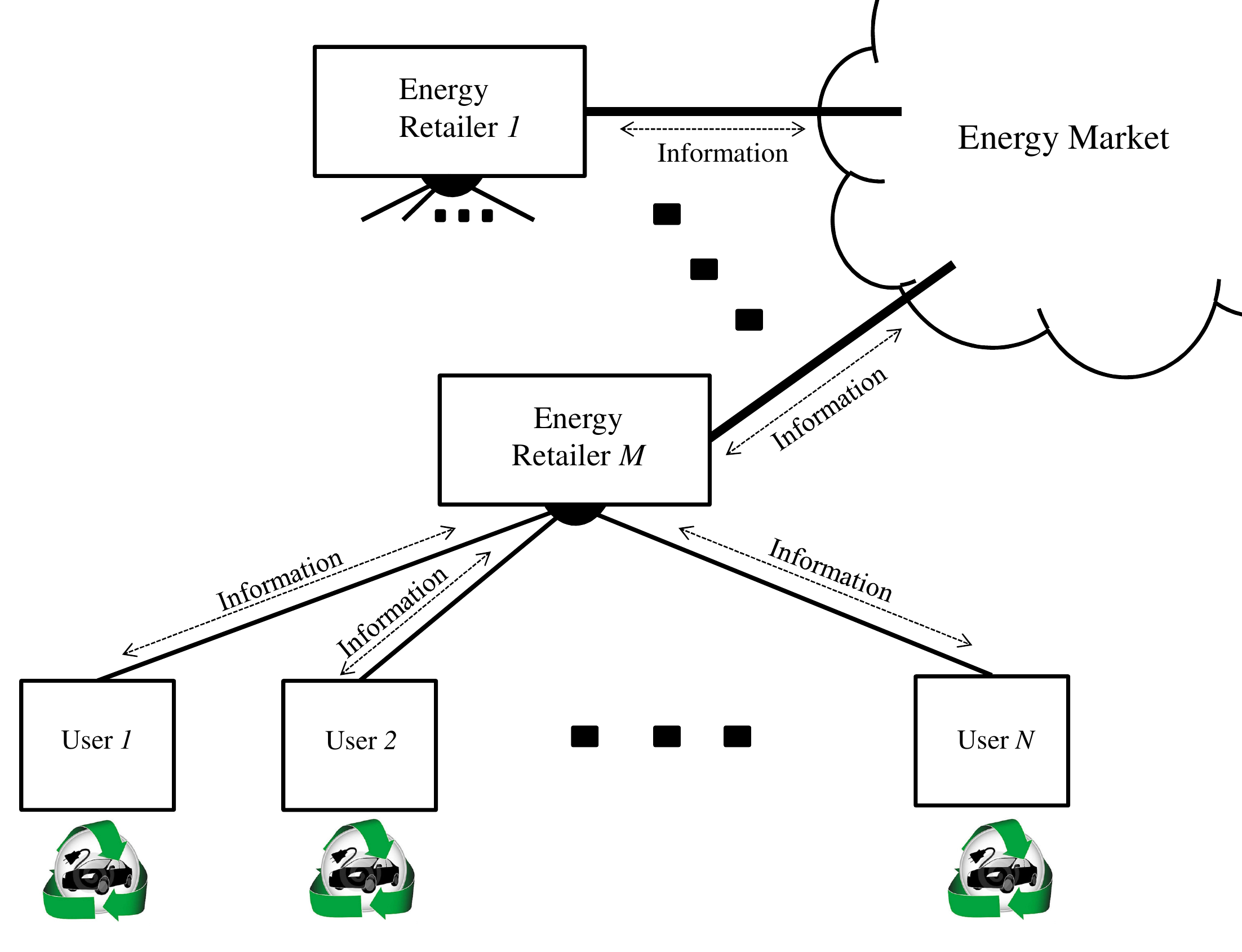} 	
	\caption{Basic model of a smart electricity system comprised of multiple users sharing one electricity retailer or an aggregator.} \vspace{-1em}
	\label{f1}
\end{figure}

Fig. \ref{f1} represents our assumed model of a smart electricity system where multiple users share one electricity retailer or an aggregator. We presume that the users' overall load consists of two distinct types of load; typical household load which normally needs \textit{on-demand} power supply, e.g. air conditioning, lighting, cooking and refrigerator, and PEV as a \textit{flexible} load.
In this model, the dotted lines show the underlying communication system while the solid lines represent the power cabling infrastructure. 
 
We assume that an electricity retailer bids to the energy market, e.g., on a DA basis. Then, based on its energy needs and the market state, it buys energy from the market at MCPs. We assume that the retailer is willing to handle its customers' PEVs' electricity assignments such that the shape of its aggregated power demand profile matches the electricity profile resulted from the successful bids in a DA market. This enables the retailer to minimize its demand from the RT market for balancing the load in the following day and accordingly reduce the overall electricity procurement cost. This cost reduction makes the energy retailer afford to offer more attractive deals to the customers in the form of pricing, rewarding, promotions, etc \cite{farst1}. 

However, since, in practice, it is not guaranteed that the shaped aggregated profile matches the retailer's purchased DA profile, the retailer often needs to reciprocate the load imbalances in the following day by referring to RT market.

On the other hand, we should notice that residential DR is desired to be carried out such that users' privacy is not violated. Therefore, DR can attract more participation from the users if it is implemented in each user's house in a decentralized fashion according to this model.   

\section{Analysis} \label{Sta}
In this section, we first formulate the electricity procurement cost for the retailer and then provide our proposed joint demand profile shaping and altering algorithm. 

In our analysis, without loss of generality, we assume a daily scheduling horizon and a time granularity of one hour for electricity markets. Let $(\textbf{l}^d,\textbf{p}^d)$ represents the pair of load $\textbf{l}^d$ and price $\textbf{p}^d$ vectors which has been cleared in the DA market. i.e.,
\begin{align}
\textbf{l}^d\triangleq [l_{1},l_{2}\dotsc,l_{H}]^T,
\end{align} 
\begin{align}
\textbf{p}^d\triangleq [p_{1},p_{2}\dotsc,p_{H}]^T,
\end{align}  
for which, the units of $l$ and $p$ are MWh and \$/MWh, respectively. Similarly, assume that $(\textbf{l}^r,\textbf{p}^r)$ represents the pair of load $\textbf{l}^r$ and price $\textbf{p}^r$ vectors which are the load imbalance and RT price vectors in the following day. The values of the elements of these vectors will be only known to the retailer at each time slot of the next day. 
Then, the overall electricity procurement cost for the next day can be formulated as follows:  
\begin{align}
C(H)=<\textbf{l}^d,\textbf{p}^d>+<\textbf{l}^r,\textbf{p}^r>,
\end{align}
in the above, $<\textbf{x},\textbf{y}>$ shows the inner product between vectors $\textbf{x}$ and $\textbf{y}$ and $C$ is the overall energy procurement cost over the scheduling horizon $H=24$.

First, given the purchased profile from the DA market by the retailer, i.e., $\textbf{l}^d$, the users individually contribute to follow this demand profile by solving the following sequential optimization problem (see Algorithm 1):
\begin{align}
\textbf{P1:} \quad \underset{\textbf{l}_{\text{PEV},n}}{\text{minimize}}
\quad <\textbf{l}_{\text{\text{PEV}},n},\textbf{l}_{\text{A},n} + \textbf{l}_{-n}-\textbf{l}^d>,
\label{P1}
\end{align}\vspace{-1.2em}
\begin{align}
&\sum_{t=\alpha_{n}}^{\beta_{n}} l^{t}_{\text{PEV},n}=E_{\text{PEV},n},\\
&|l^{t}_{\text{PEV},n}| \leq p_{max},\quad \forall t\in\mathbb{T}^P_{\text{PEV},n}, \\ 
& l^{t}_{\text{PEV},n}=0,  \forall t\notin\mathbb{T}^P_{\text{PEV},n},\\
& SOC^{t=\alpha}_{\text{PEV},n}+\sum_{k=\alpha+1}^{t}l'^k_{\text{PEV},n}\geq 0.2 \times C_{\text{PEV},n}, 
\forall t\in\mathbb{T}^P_{\text{PEV},n}.
\end{align}
\noindent Here, $\textbf{l}_{\text{PEV},n}$ and $\textbf{l}_{\text{A},n}$ show the energy assignment vector for user n's PEV and the aggregated load from its household appliances, respectively. $E_{\text{PEV},n}$ is the $n^{\text{th}}$ user's required energy to be allocated to its PEV which is associated with the total required charging time $T_{\text{PEV},n}$ as follows: 
\begin{align}
E_{\text{PEV},n}=a\times T_{\text{PEV},n},
\end{align} 
\noindent where $a$ is the charging power rate. Likewise, $\alpha_{n}$ and $\beta_{n}$ represent the arrival time and departure time of the PEV. Furthermore, $|l^{t}_{\text{PEV},n}| \leq p_{max}$ limits the maximum power that can be delivered to/from the PEV and $\mathbb{T}^P_{\text{PEV},n}$ represents the permissible charging time set or simply the set of time slots during the PEV's \textit{connection time} to the power grid. Additionally, $\textbf{l}_{-n}$ is the aggregated power profile from other $N-1$ users described as follows: 
\begin{align}
\textbf{l}_{-n}=\sum_{\substack{i \in {N} \\  i \neq n}}  (\textbf{l}_{\text{PEV},i}+\textbf{l}_{\text{A},i}). 
\label{others}
\end{align}
In (8), $C_{\text{PEV},n}$ is the total storage capacity of the user n's PEV and we assumed that in case of employing V2G in the system, PEV's $SOC$ should not fall below 20\% of that total capacity in order to make sure that the adverse impacts on PEV's battery due to complete depletion would not take place.   
 
Second, knowing the fact that $(\textbf{l}^r,\textbf{p}^r)$ are unknown to the retailer a priori, at each time slot $t_0$ of the next day after getting this information, the retailer may decide to reshape the previous shaped demand profile to minimize its RT electricity purchase to balance the load and even sell back some of its pre-purchased electricity from DA market if the price rises significantly due to the state of the RT market or contingencies that may increase the price unexpectedly at specific time slots of the day. 

\begin{figure}[t]
	\centering
	\includegraphics[width=\columnwidth,height=1.7in]{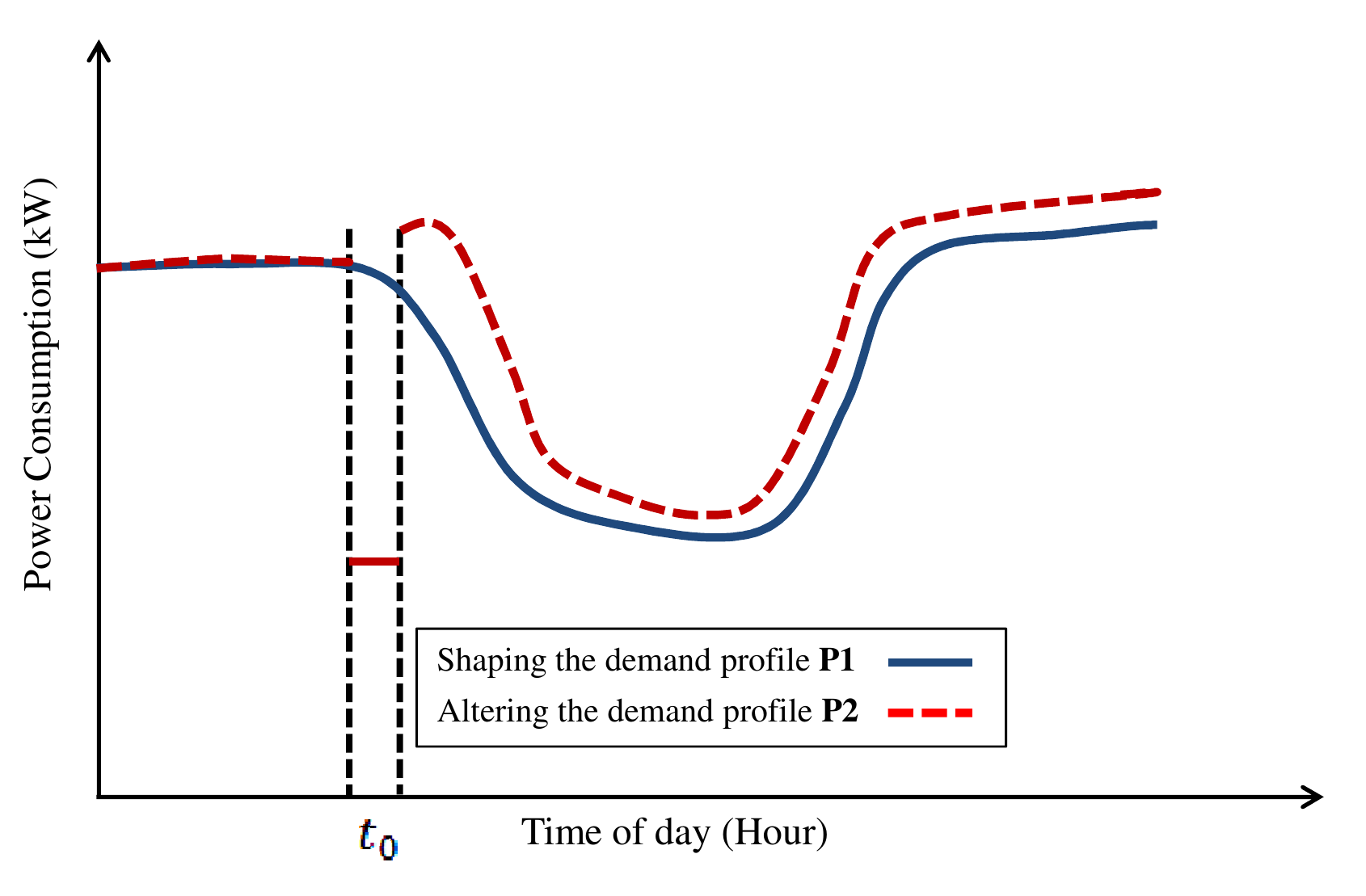}
	\caption{Reshaping the demand profile upon a contingency.} 
	\label{f4}
\end{figure}

\begin{figure}[t] 
	\centering
	\includegraphics[width=\columnwidth,height=1.7in]{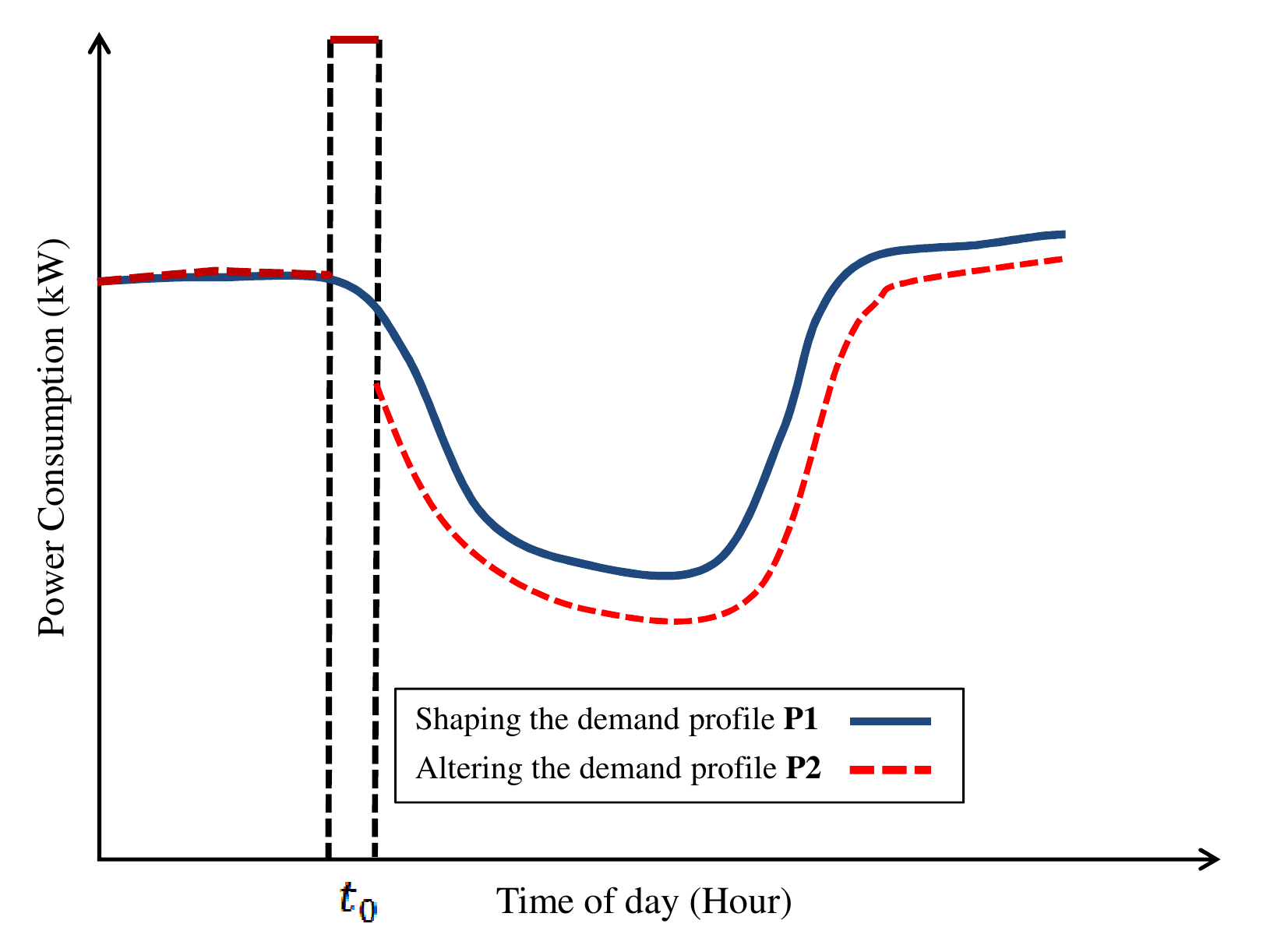}
	\caption{Reshaping the demand profile upon a contingency.} \vspace{-1em}
	\label{f5}
\end{figure}

As the chance of staying the price that much during the next few hours of the day is low \cite{poterba1988mean}. Reshaping the load profile by lowering the electricity consumption at that time slot and purchasing electricity at the next time slots will yield to a lower electricity procurement cost. This is also true for purchasing electricity at those time slots that price, unexpectedly, falls down significantly. The retailer may buy extra electricity at those specific time slots (based on its overall storage capacity coming from connected PEVs) and reshape the demand for the next hours, c.f., Fig. \ref{f4} and Fig. \ref{f5}.

We should note that the retailer is assumed to be allowed to use the \textit{existing} flexibility (from each PEV) and the diversity (coming from all the users). In other words, in our case the electricity consumption behaviours of the users (their PEV usage patterns) are not to be changed.



Then, the following linear multi objective programming (MOP) allows demand altering along with pursuing the shape of the pre-purchased electricity from DA market (see Algorithm 1): 
\begin{align}
\textbf{P2:}  \quad \underset{\textbf{l}'_{PEV,n}}{\text{min}} \quad
(\lambda <\textbf{l}'_{PEV,n},\textbf{l}_{A,n} + \textbf{l}'_{-n}-\textbf{l}^*_{N}>\\+ (1-\lambda) (l'^{t_0}_{-n}+l^{t_0}_{A,n}+l'^{t_0}_{PEV,n})),
\label{P2}
\end{align}\vspace{-1.7em}
\begin{align}
&[l^{'1}_{PEV,n},\cdots,l^{'t_0-1}_{\text{PEV},n}]=[l^{1}_{\text{PEV},n},\cdots,l^{t_0-1}_{\text{PEV},n}]   ,\\
&\sum_{t=t_0}^{\beta_{n}} l'^{t}_{\text{PEV},n}=E_{\text{PEV},n}-\sum_{t=\alpha_{n}}^{t_0-1} l^{t}_{\text{PEV},n},  \\ 
&|l'^{t}_{\text{PEV},n}| \leq p_{max},\quad \forall t\in\mathbb{T}^P_{\text{PEV},n}, \\
& l'^{t}_{\text{PEV},n}=0, \quad \forall t\notin \mathbb{T}^P_{\text{PEV},n},\\ 
& SOC^{t=t_0-1}_{\text{PEV},n}+\sum_{k=t_0}^{t}l'^k_{\text{PEV},n}\geq 0.2 \times C_{\text{PEV},n}, 
\forall t\in\mathbb{T}^P_{\text{PEV},n},
\end{align}
where $\lambda$ and $(1-\lambda)$ are the weights of the objective functions. Furthermore, $t_0$ is the request time slot for \textit{altering the demand}. In (17), $SOC^{t=t_0-1}$ is the $SOC$ of the user n's PEV at time slot $t_0$ and, similar to \textbf{P1}, we assumed that in case of employing V2G in the system, the $SOC$ should not fall below 20\% of the total PEV's storage capacity.   

The convergence criterion in Algorithm 1 can be simply assumed as a desired number of iterations of updating all users' demand profiles or it can be determined to be derived from measuring mean square error (MSE) between two subsequent iterations upon achieving aggregated demand profiles. As we have discussed in \cite{farst1}, a convergence is guaranteed to be obtained. Furthermore, users' contribution can be modelled as a cooperative game with complete information wherein a Nash equilibrium exists \cite{farst1}.      


\begin{algorithm}[t]
	\begin{algorithmic}[1]
		\caption{Joint Demand Shaping \& Altering}\label{demshape}
		\State All $N$ users initialize their respective load profiles over the scheduling horizon based on their respective demands, i.e., $\textbf{l}_{n}$ for $n=1,\dotsc, N$.
		\State All $N$ users send their initialized load profiles to the retailer.
		\While {\textit{not reaching convergence}}
		\For {$n=1$ to $N$}
		\State  The retailer calculates the state information $\textbf{l}_{-n}$ according to (\ref{others}) for user $n$.
		\State The retailer sends $(\textbf{l}_{-n}-\textbf{l}^d)$ to user $n$.
		\State User $n$ solves problem \textbf{P1} and updates its load profile $\textbf{l}_{n}$.
		\State User $n$ sends back the new demand profile to the retailer. 
		\State The retailer updates $\textbf{l}_{n}$.
		\EndFor 
		\State \textbf{end for}
		\EndWhile 
		\State \textbf{end while}
		\State The retailer receives information from RT market. 
		\State The retailer finds $t_0$.
		\While {\textit{not reaching convergence}}
		\For {$n=1$ to $N$}
		\State The retailer sends demand altering signal at time slot $t_0$ to user $n$.
		\State User $n$ solves problem \textbf{P2} and updates its load profile $\textbf{l}^{'}_{n}$.
		\State User $n$ sends back the new demand profile to the retailer. 
		\State The retailer updates $\textbf{l}^{'}_{n}$.
		\State  The retailer calculates the state information $\textbf{l}^{'}_{-n}$ according to (\ref{others}) for user $n$.
		\EndFor 
		\State \textbf{end for}
		\EndWhile 
		\State \textbf{end while}
	\end{algorithmic}
\end{algorithm}

\section{Simulation Results} \label{SR}
In this section, we evaluate the effectiveness of the DR techniques described in the previous section, through extensive computer simulations. In our simulations, the number of users, \textbf{N}, is 1,000 and the optimization horizon is considered to be a day, i.e., 24 hours for a DA programming scenario and a time granularity of one hour for each time slot. 

For the PEVs usage patterns, our data and distributions are based on 2009 NHTS data \cite{NHTS}. We considered new standard outlets, NEMA 5-15, providing 1.8 kW power. Furthermore, $SOC$ for each PEV at the arrival time is as follows in percentage points: 
\begin{equation}
SOC^{t=\alpha}_{\text{PEV},n}=100 \times (1-\frac{E_{\text{PEV},n}}{24}).
\end{equation} 
In other words, we assumed that PEVs are needed to be fully charged by their respective next departure time. Additionally, we considered 24 kWh energy storage capacity for PEVs according to Nissan Leaf model \cite{leaf}. We adopted the PJM interconnection electricity market pricing data for both DA and RT markets in 2014 \cite{PJM}.

\begin{figure} 
	\centering
		\includegraphics[width=\columnwidth]{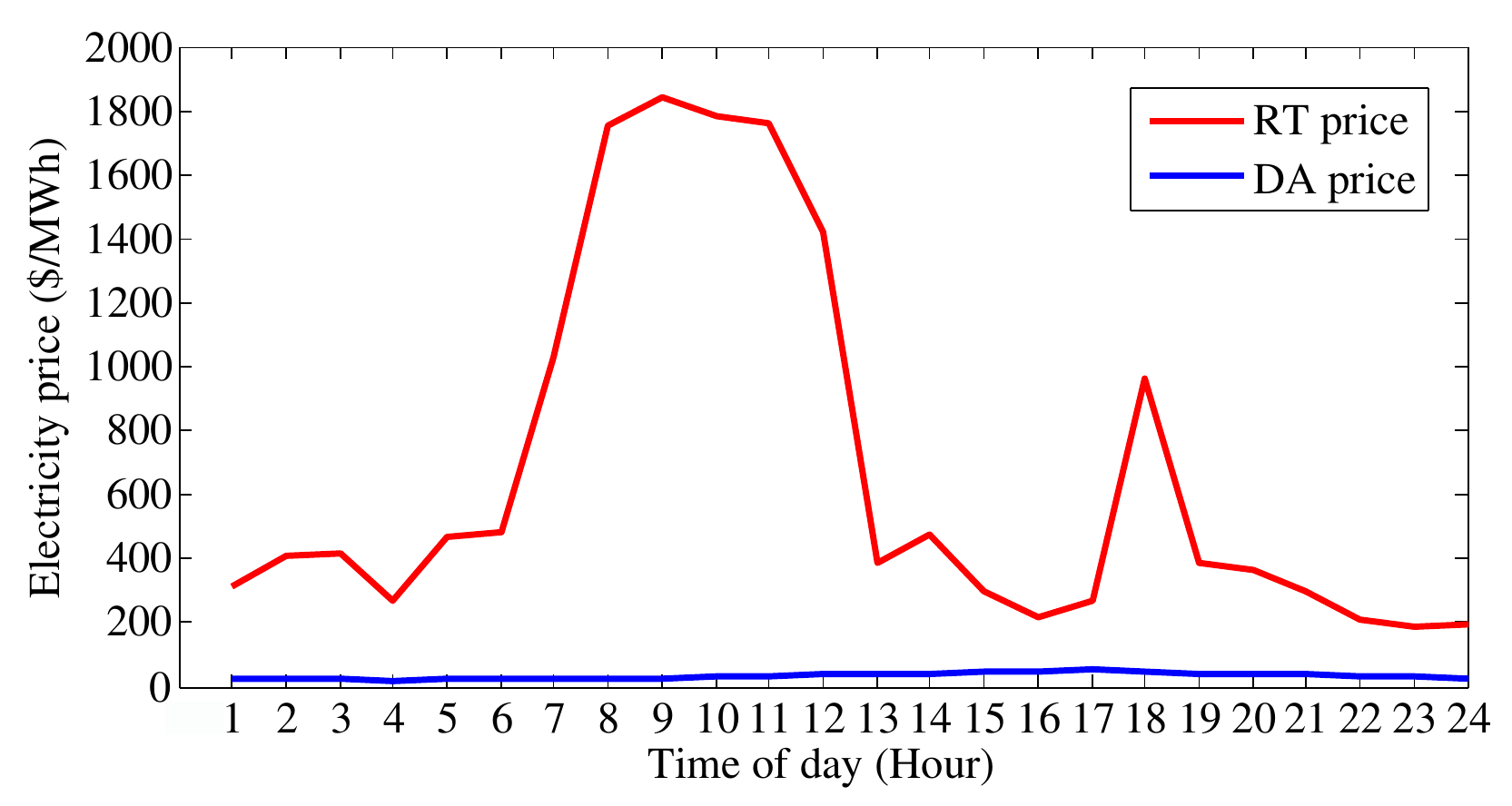}
	\caption[]{DA and RT prices for March 9 2014 in PJM interconnection electricity market (the number of days with such behaviour in 2014 is quite considerable).} 
	\label{price}
\end{figure}

Fig. \ref{price} shows the DA and RT prices for March 9 2014 as an example. We chose this day since it had the highest peak in RT prices throughout the year in PJM and can be assumed as what we referred to as a black swan event in the market. As it can be observed, the RT price has a substantial peak at 9 A.M. around which the price is still unexpectedly high for 5 hours. The retailer can choose to proceed for demand altering program \textbf{P2} upon receiving this pricing information. 

\begin{figure} 
	\centering
	\includegraphics[width=\columnwidth]{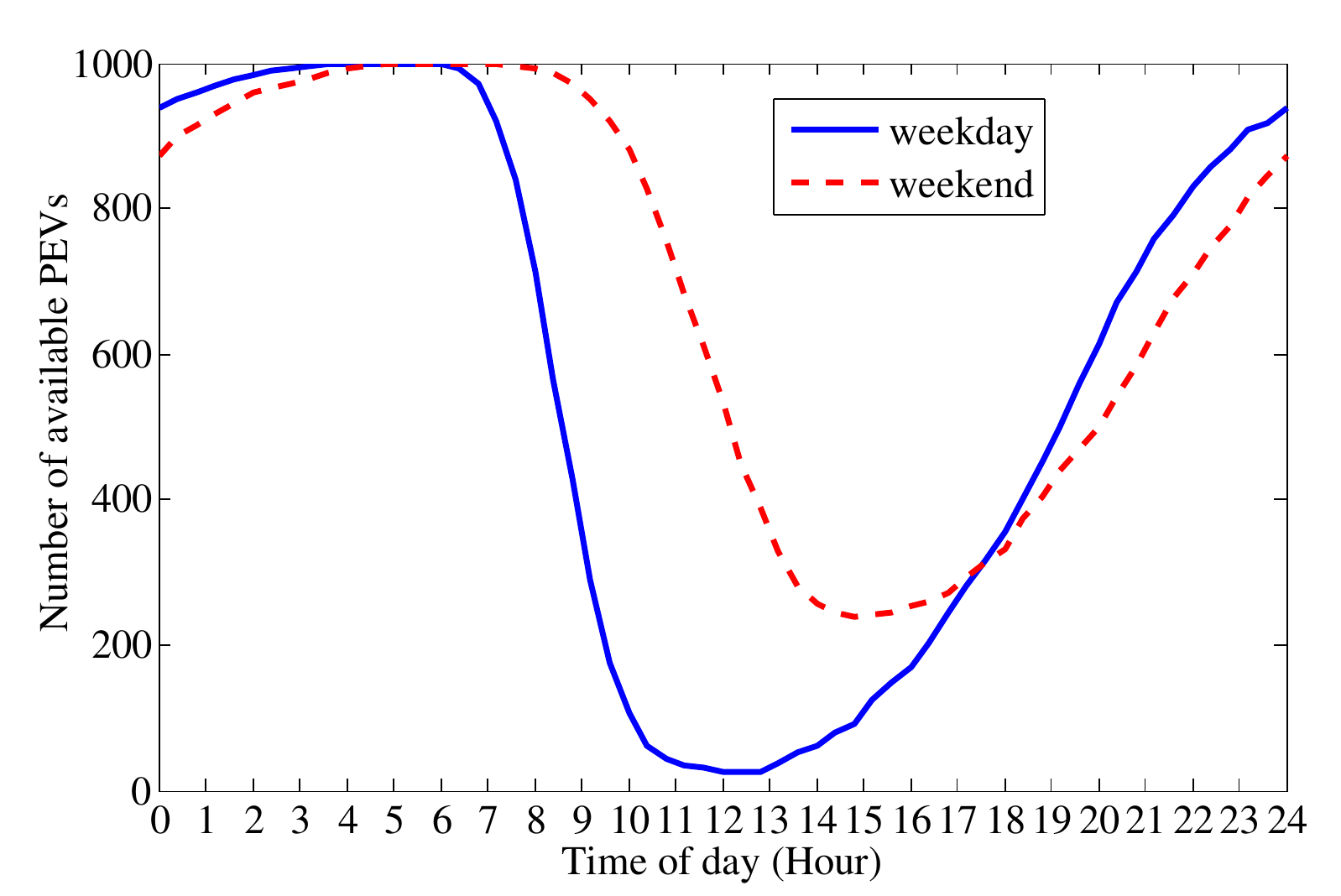}
	\caption[]{The number of PEVs available at users' dwellings for different hours of a day according to NHTS \cite{NHTS}.}
	 
	\label{ava}
\end{figure}

The amount of success in demand altering depends on the availability of PEVs at users' dwellings which is dissimilar for different hours of a day and different days of a week, see Fig. \ref{ava}.

\begin{figure} 
	\centering
	\includegraphics[width=\columnwidth,height=3.3in]{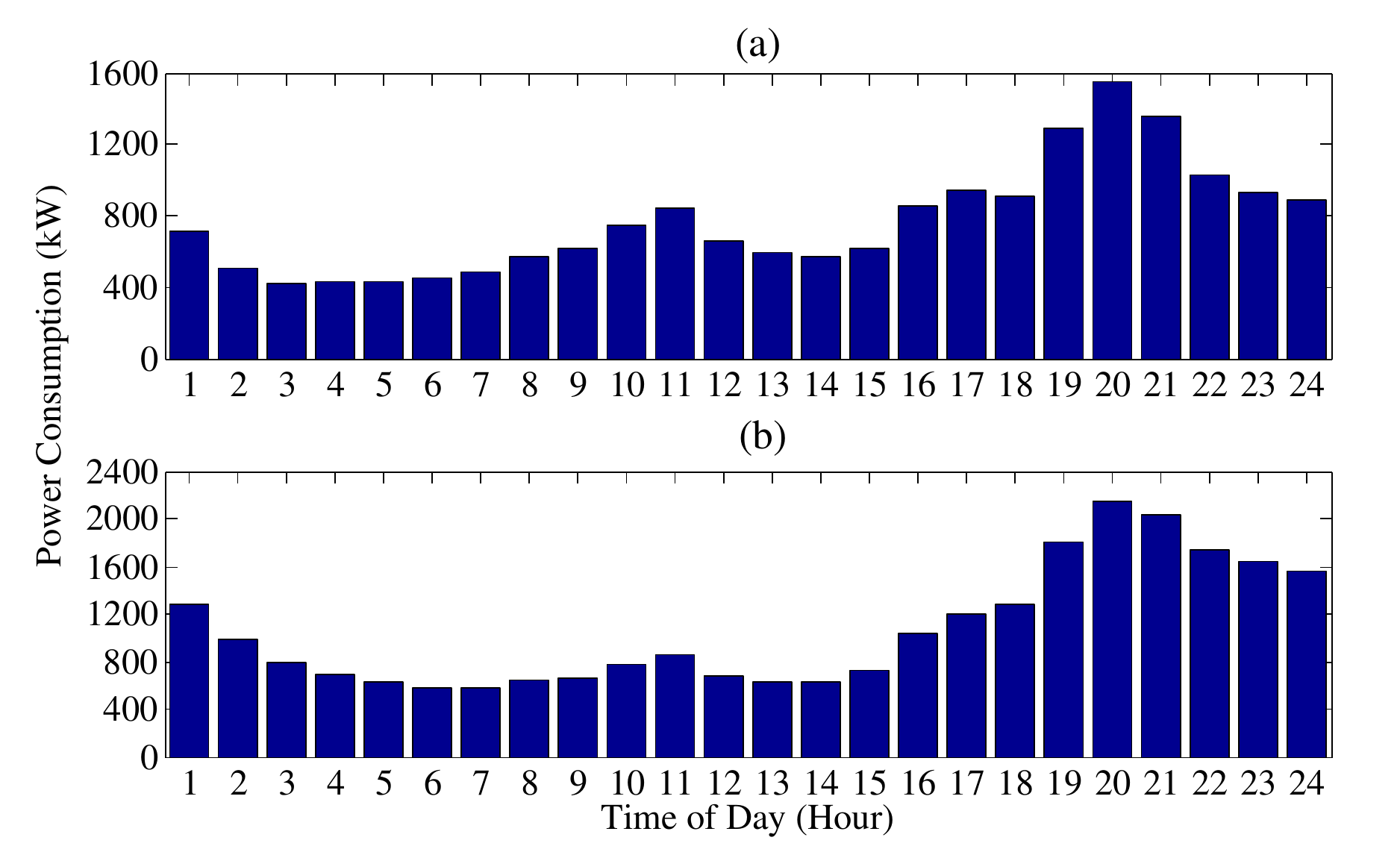}
	\caption[]{Electricity demand profile from (a) only normal household appliances, i.e., without PEVs and (b) the overall electricity demand profile when users use PEVs with different usage patterns based on NHTS data.} 
	\label{infev}
\end{figure}

Next, we examine the DR scheme introduced in Algorithm 1. Fig. \ref{infev} shows the electricity demand profile from only typical household appliances, i.e., without PEVs and the overall electricity demand profile when users use PEVs with different usage patterns based on NHTS data.

\begin{figure} 
	\centering
	\includegraphics[width=\columnwidth]{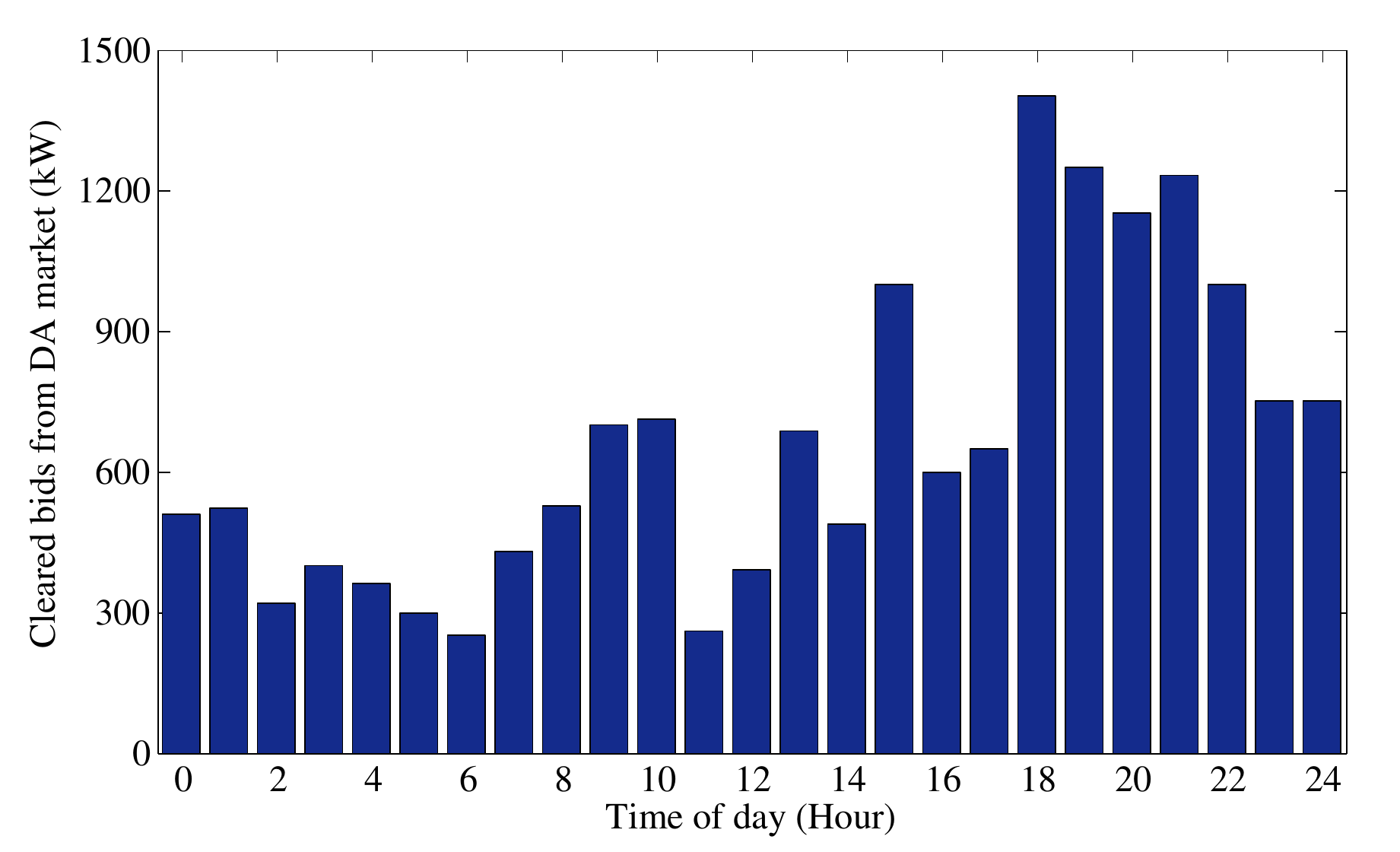}
	\caption[]{The assumed electricity profile purchased by the retailer from the DA market by the bids that could be cleared.} 
	\label{dap}
\end{figure}

\begin{figure} 
	\centering
	\includegraphics[width=\columnwidth,height=3.3in]{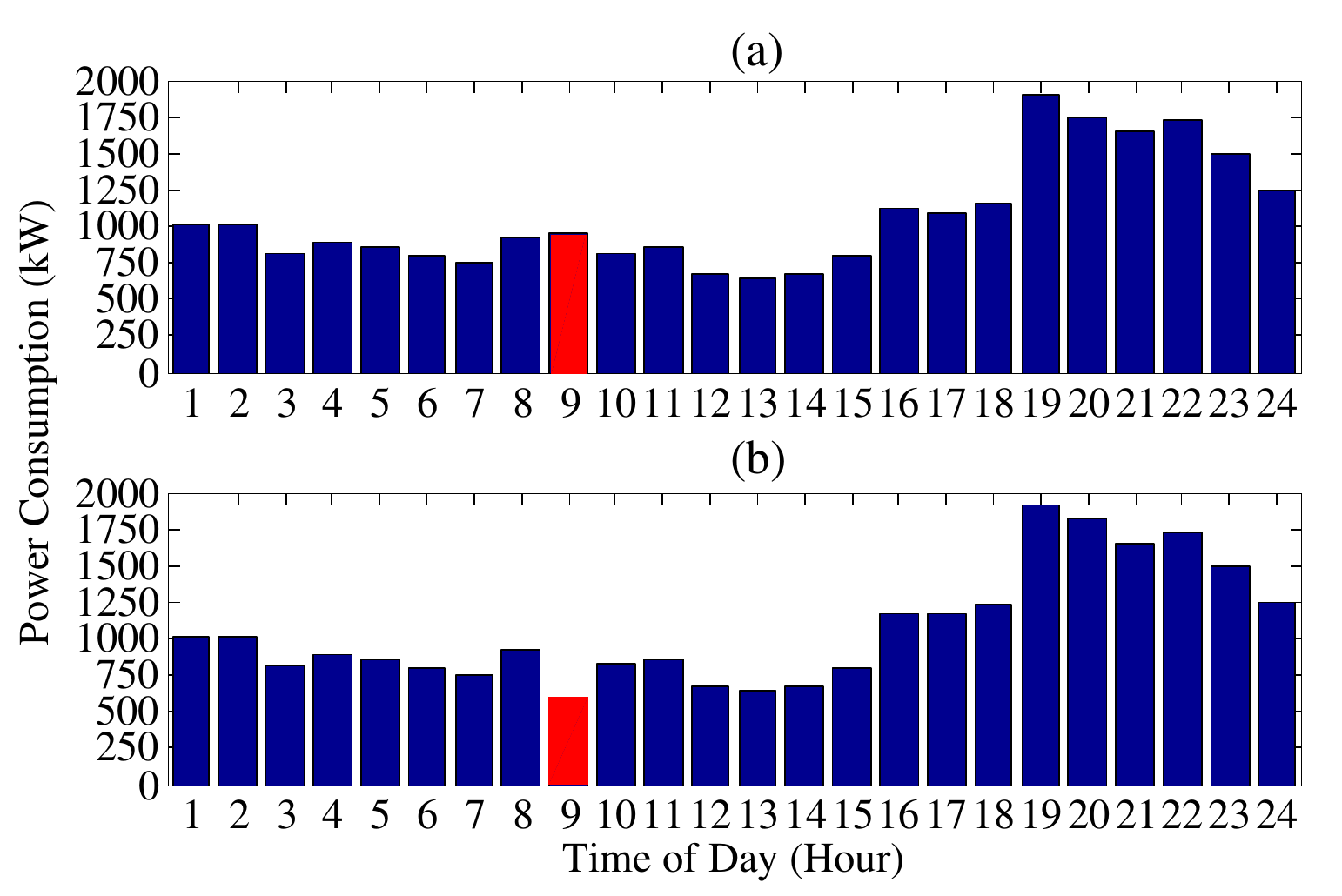}
	\caption[]{Load profiles after using Algorithm 1 for (a) shaping and then (b) altering the demand at 9 A.M.} 
	\label{shre}
\end{figure}

We assume that Fig. \ref{dap} shows the electricity profile purchased by the retailer from the DA market. In other words, it shows the bids that could be cleared in the market. The results of Algorithm 1 for shaping the demand profile and joint shaping and altering the demand is depicted in Fig. \ref{shre}. For this, we assumed $\lambda=0.5$ in \ref{P2}. It can be observed that at 9 A.M. when the highest RT price occurs, demand altering can reduces the aggregated demand from 948.2 kWh to 599.3 kWh, i.e., we can obtain almost \%37 reduction in the overall demand. 

In our simulations, convergence has been attained only after one single iteration of updating \textit{all} users' electricity demand profiles for both \textbf{P1} and \textbf{P2} in Algorithm 1.  

It should be emphasized that this could be achieved since at that hour of the day we had almost 330 V2G enabled PEVs available at users' dwellings. Different results would be obtained for the other hours of that day. Especially, it is obvious that the amount of cost savings would be dissimilar for a weekday and a weekend day, see Fig. \ref{ava}.     

\begin{table}
	
	\renewcommand{\arraystretch}{2.5}
	\caption{\textlcsc{Overall energy procurement costs for the retailer}}
	\label{tab}
	\centering
	\resizebox{\columnwidth}{!}{  
		\begin{tabular}{|c|c|c|c|c|}
			\hline
			\normalsize{Case} & \normalsize{Overall ideal cost (\$)} & 
			\normalsize{Overall real cost (\$)}
			&
			\normalsize{Overall cost after \textbf{P1}} (\$)
			& \normalsize{Overall cost after \textbf{P2} (\$)} \\ 
			\hline
			\normalsize{1} & \normalsize{\$674.4} &  \normalsize{\$2808.9} & \normalsize{N/A} & \normalsize{N/A}\\
			\hline
			\normalsize{2}  &  \normalsize{\$920.4} & \normalsize{\$5865.1} & \normalsize{\$4775.6} & \normalsize{\$4308.9}\\
			\hline
		\end{tabular}\vspace{-2.5 em}
	}
	\label{tab1}
\end{table}

In Table \ref{tab1}, we compare the overall energy procurement costs for the retailer for two cases: \textit{case 1)} when there is no PEV in the system and \textit{case 2)} when all the users possess PEVs with their respective usage patterns. 

It can be noticed that in case 1, if the retailer could be absolutely successful in bidding to the DA market, i.e., there would not be any need to purchase electricity from the RT market, the overall cost is so much lower for our assumed pricing in Fig. \ref{price}. In a more realistic case, when the retailer's bidding to the DA market is assumed to be according Fig. \ref{dap}, and the retailer is required to balance the load, the overall cost is much higher. 

In this case, we also assumed that the retailer can sell back its extra load purchased earlier from the DA market to the RT market at the same price offered by the RT market. This happens a few times such as at 4 P.M. Obviously, demand shaping and demand altering in this case is not applicable (N/A) as there is no PEV in the system. 

For the second case, when users possess PEVs, for the ideal bidding the overall cost increases by almost \%37 to supply electricity to the PEVs whereas for the realistic bidding it becomes more than double. 

When demand shaping \textbf{P1} is employed, this overall cost reduces by around \%18. Furthermore, when joint demand shaping and altering in Algorithm 1 is used, it decreases further by almost \%10.


\section{Conclusion and Future Work} \label{Con}

In this paper, we proposed a fast converging and decentralized joint demand shaping and altering algorithm for managing V2G enabled PEVs' electricity assignments (charging and discharging) to reduce the overall electricity procurement cost for an electricity retailer. Our proposed algorithm uses demand shaping and demand altering for the DA and the RT markets, respectively. In particular, when the power system has high penetration of intermittent energy resources, demand altering is crucial due to likely contingencies and hence more RT price fluctuations. In our simulations' results, we considered a specific day in 2014 which had high unexpected RT prices for some hours in PJM interconnection electricity market. We showed that significant overall cost savings (up to \%28) for a retailer bidding to this electricity market can be achieved by using our proposed algorithm. This allows the retailer to offer better deals to the customers and expand its market capacity. Therefore, customers can enjoy lower electricity bills as well.  

The work presented in this paper can be extended in various ways. First, the overall amount of flexibility offered by PEVs for demand altering throughout different hours of a day can be obtained. Second, the retailer can determine a threshold based on which it proceeds for demand altering. This can be done via stochastic dynamic programming (SDP) and data mining. Third, recently, PEVs' battery storage capacity has been expanded up to 60-85 kWh, e.g., for Tesla model S and model X \cite{Tesla}, which can provide much more elasticity for demand shaping and demand altering and hence reduce the electricity costs further. Finally, demand altering for the RT market does not necessarily need to take place at only one specific hour of the day as the high unexpected price can remain large for a few hours. We are addressing all these points together with some more detailed analyses in the follow up journal paper of this work.     

%
%

\bibliographystyle{IEEEtran} 
\bibliography{IEEEabrv,myBIB}

\end{document}